\documentclass[11pt]{amsart}
\usepackage[T1]{fontenc}
\usepackage[utf8]{inputenc}
\usepackage{lmodern}
\usepackage{amsmath,amssymb,amsthm}
\usepackage[margin=1in]{geometry}
\usepackage{microtype}
\usepackage[hidelinks]{hyperref}
\newenvironment{namedtheorem}[1]{\par\addvspace{8pt}\noindent\textbf{#1.}\ \itshape}{\par\addvspace{5pt}}
\newenvironment{nameddefinition}[1]{\par\addvspace{8pt}\noindent\textbf{#1.}\ \normalfont}{\par\addvspace{5pt}}
\newenvironment{namedremark}[1]{\par\addvspace{8pt}\noindent\textit{#1.}\ \normalfont}{\par\addvspace{5pt}}
\title[Profile of a representation variety]{Algebro-Geometric Invariants of Finitely Generated Groups\\(The Profile of a Representation Variety)}
\author{S. Liriano}
\email{SAL21458@yahoo.com}
\author{S. Majewicz}
\email{smajewicz@kbcc.cuny.edu}
\date{Corrected working version, second audit, October 8, 2026}
\thanks{We thank Ray Hoobler for his many useful comments. Also, thanks to David Smith of the NYPL HSS Division.}
\begin{document}

\begin{abstract}
If $G$ is a finitely generated group (fg), and $A$ is a complex affine algebraic
 group, then the space
 $Hom(G,A)$ admits the structure of an algebraic variety, not necessarily irreducible,
 denoted briefly by $R_A(G)$. The invariants of $R_A(G)$  are
 independent of the
 finite set $X$ of generators for $G$ chosen. In this communication
 we define the profile function, $P_d(R_A(G))$, of the representation
 variety over an algebraic group $A$ of a fg group $G$ to be
 $P_d(R_A(G))=(N_d(R_A(G)),\dots,N_0(R_A(G)))$, where $N_i(R_A(G))$ stands for the
 number of  irreducible components  of $R_A(G)$ of dimension
 $i$, $0 \leq i \leq d$, and $d=Dim(R_A(G))$. We apply this invariant to
 the study of fg groups and prove various results.
 In particular, we show that if $G$ is the fundamental
 group of a closed orientable surface of genus $g \ge 1$,  then
 
\begin{equation*}
P_d(R_{SL(2,\mathbb C)}(G))\neq P_d(R_{PSL(2,\mathbb C)}(G)).\tag {**}
\end{equation*}

 In (**) the index $d$ is a common upper bound for the two dimensions.
 We also show that (**) holds
 for $G$ a torus knot group with presentation $<x,y;x^p=y^t>$ where
 $p,t$ are both $\geq 3$, and that (**) also holds
when $G$ is the
 fundamental group of a closed non-orientable surface of genus $g \ge 3$.
 Additionally, we show that if a group $G$ can be  $n+1$ generated and admits a
 presentation $<x_1,\dots,x_n,y \, ; \, W=y^p>$, where $p>1$, $W$ is a
 non-trivial word in $F_n=<x_1,\dots,x_n>$, and $A=PSL(2,\mathbb C),$ then
 $Dim(R_{A}(G))=Max\{3n,\,Dim(R_{A}(G'))+2 \} \le 3n+1$,
 where $G'=<x_1,\dots,x_n \, ; \, W=1>$. We also give a condition
 guaranteeing that the algebraic variety is reducible.
\medskip

\end{abstract}
\maketitle

\begin{center}
\small\textbf{Version note.} This is a corrected working version prepared from
the source of arXiv:0804.0460v1 (2008), retaining its theorem numbering.
The paper was subsequently published under the title
\emph{Algebro-Geometric Invariants of Groups (The Dimension Sequence of a
Representation Variety)}, \emph{International Journal of Algebra and
Computation} \textbf{21} (2011), 595--614. This working version has not
been collated against the journal text and is not a journal-issued erratum.
New extensions are kept in the separate revision notes.
\end{center}
\medskip

 \section*{Historical Remark}

 \medskip

 We mention briefly only a few
 developments bearing some
 resemblance to this work, keeping in mind
  that this account
  is anything but thorough; in passing, let us remark that  the study
  of spaces of representations of fg groups can be traced back to
 the work of H. Poincare \cite {LM}. We begin  with the work
 of  W. Goldman, who in
 the eighties produced results
giving the number of connected components of $Hom(G, L)$, where
$L$ is a Lie group and $G$ an orientable surface group; for
example, in \cite {GW} a result was laid out yielding the number
of connected components in the space of representations of the
fundamental group of an oriented surface in the $n$-fold covering
group of $PSL(2,\mathbb R)$ in terms of $n$ and the genus of the
surface. We recall, however, that counting connected
 components of a complex variety is fundamentally
 different from counting its irreducible components. For example, consider the vanishing
 set $V$ of  $xy=0$ in $\mathbb C^2$. It is easy to see that
$V$ has two   irreducible components and only one connected
component.
  More recently, in the nineties,
  Rapinchuk, et al, \cite{RBC} produced a paper
  where amongst numerous interesting results, they established
  the absolute irreducibility of
  $R_{\Lambda}(\Gamma_g)$, the
  representation variety of the fundamental group of a closed
  orientable surface of genus $g$ for $\Lambda $ either $GL(n,K)$, or
   $SL(n,K)$, where $K$ is a field of characteristic $0$; they also gave
   formulas for the  dimension of their corresponding
  representation varieties  involving the genus, and the dimension
  of the relevant algebraic group.
  Subsequently, Benyash-Krivetz and Chernousov, \cite {BC},
  obtained formulas yielding
the number of
 irreducible components of $R_{\Lambda}(\Gamma_g)$, where $\Gamma_g$ is a
 closed
 non-orientable surface group of genus $g$, and $\Lambda$ is either
 $GL(n,K)$, or $SL(n,K)$, where $K$ is an algebraically closed field
 of characteristic $0$. At around the same time,
 Liriano in \cite{L3} produced a formula giving the number
 of four
 dimensional irreducible components of  $R_{SL(2,\mathbb C)}(G)$ for groups
 $G$ having presentation $<x,y \, ; \, x^p=y^t>$ where \, $p,\,t$ are positive
  integers, and pointed
out that when
 such a group
 is the fundamental group of a torus knot complement in $\mathbb R^3$,
 that its number of four dimensional irreducible components equals the
 genus of the corresponding torus knot.

 \section*{Introduction}\enlargethispage{2pt}

\medskip
 The study of invariants of representation varieties over $SL(2,\mathbb C)$
 of parafree groups \footnote {A group $G$ is parafree of rank $r$ if it is residually nilpotent
 and has the same lower central sequence as a free group of rank $r$.
 Its deviation is the minimum number of generators of $G$ minus $r$. For a more detailed
 account on parafree groups consult \cite {B1}, \cite {B2}.}
 has provided ample evidence to justify the
 claim that invariants of groups associated with spaces of
 representation over algebraic groups can
 be
 highly sensitive, even when launched in classes of groups exhibiting
 pronounced structural  affinities \cite {L4}, \cite {L2}. In this
 paper an invariant, $P_n(R_A(G))$, is introduced and we call it {\it
 \lq\lq the profile function of the
 representation
 variety of $G$ over the algebraic group $A$"}, or more swiftly,
  {\it \lq\lq the profile of $G$ over $A$\rq\rq }. This invariant
  that counts all  irreducible
 components of all dimensions in $R_A(G)$ is then deployed, in what may be described as
 \lq doubly
 close settings\rq;
 by that it is meant that not
 only is the invariant studied over groups that are known to be
 structurally very close to each
other, for example, free groups and  surface groups, (or more
generally, the class of fully residually free groups), but also
over algebraic groups that are also structurally similar; for
example, $SL(2,\mathbb C)$ and $PSL(2,\mathbb C)$. Even when we observe
this doubly close condition, if you will,  rather remarkable
divergence in the respective values of the invariant is obtained.
\medskip

 More formally, define the profile function $P_d(R_A(G))$ of the
 representation variety over an algebraic group $A$ for a fg
 group $G$  to be
 
\begin{equation*}
P_d(R_A(G))=(N_d(R_A(G)),\dots,N_0(R_A(G))),
\end{equation*}
 where $N_i(R_A(G))$ stands
 for the
 number of  irreducible components  of $R_A(G)$ of dimension
 $i$, $0 \leq i \leq d$, and $d=Dim(R_A(G))$. From the start, let's
 make it clear that it is possible
 to take the profile function $P_n(R_A(G))$ at any integral value
 $n\ge Dim(R_A(G))$, padding with zeros on the left. When two
 varieties are compared, we use a common index at least as large as both
 dimensions. An unindexed profile may instead be written with its leading
 zeros omitted. In this paper
 the following results are
 obtained:

 \begin{namedtheorem}{Theorem A} Let $w$ be a non-trivial freely reduced
 word in the
 commutator subgroup of $F_n=<x_1,\dots,x_n>$, and let
  $G=<x_1,\dots,x_n \, ; \, w=1>$. If $R_{SL(2,\mathbb C)}(G)$ is an
  irreducible variety and
$V_{-1}=\{\rho \, \vert \, \rho \in R_{SL(2,\mathbb
C)}(F_n),\rho(w)=-I\}\neq \emptyset$, then $P_d(R_{SL(2,\mathbb
C)}(G))\neq P_d(R_{PSL(2,\mathbb C)}(G)).$
\end{namedtheorem}

\begin{namedtheorem}{Theorem B} Let $w$ be a non-trivial freely reduced word
in the free group on $\{x_1,\dots,x_n\}$ with even exponent sum on
each generator, and with exponent sum not equal to zero on at
least one generator. Suppose
  $G=<x_1,\dots,x_n \, ; \, w=1>$. If $R_{SL(2,\mathbb C)}(G)$ is an
  irreducible variety, then
  $P_d(R_{SL(2,\mathbb C)}(G))\neq P_d(R_{PSL(2, \mathbb C)}(G)).$
\end{namedtheorem}

An immediate consequence of the above results, together with
results found in \cite{RBC} and \cite {BC}, is that when $G$ is
the fundamental group of a closed orientable, or a closed
non-orientable surface of genus $g\geq 1$,\, $g\geq 3$,
respectively, then $R_{PSL(2,\mathbb C)}(G)$ is reducible and
$P_d(R_{SL(2,\mathbb C)}(G))\neq P_d(R_{PSL(2,\mathbb C)}(G))$. In \cite
{BC} it was shown that $R_{SL(2,\mathbb C)}(G)$ is irreducible when
$G$ is a closed non-orientable surface group of genus $g\geq 3$,
and in \cite{RBC} it was shown that $R_{SL(2,\mathbb C)}(G)$ is
irreducible when $G$ is a closed orientable surface group of genus
$g\geq 1$.
\medskip
Furthermore, we show that if $G$ is a torus knot group with
presentation $<a,b \, ; \, a^p=b^t>$, where both $p,t \geq 3,$
then $P_d(R_{SL(2,\mathbb C)}(G))\neq P_d(R_{PSL(2,\mathbb C)}(G))$. In
particular,
 the genus of these torus knots does not equal
 $N_4(R_{PSL(2,\mathbb C)}(G))$, whereas \cite{L3} identifies the genus with
 $N_4(R_{SL(2,\mathbb C)}(G))$.

\medskip
To illustrate the  sudden changes of invariants that can arise
when looking at the representation variety of a group $G$ over one
algebraic group as opposed to another, consider the following: let
$G=<x,y \, ; \, [x,y^2]y^2=1>$, and define $[a,b]=aba^{-1}b^{-1}$.
Then $Dim(R_{SL(2,\mathbb C)}(G))=3$, and $P_3(R_{SL(2,\mathbb
C)}(G))=(2,0,0,0)$, while $Dim(R_{PSL(2,\mathbb C)}(G))=5$, and
$P_5(R_{PSL(2,\mathbb C)}(G))=(1,0,1,0,0,0).$
With the stated commutator convention the relator simplifies to
$xy^2x^{-1}$, so $G\cong\mathbb Z*\mathbb Z_2$. The profiles follow
from the cyclic-group calculations and the free-product identity below.
More generally, let $G_r$ be the free product of $r$ copies of
$\mathbb Z_2$. Then
\[
R_{SL(2,\mathbb C)}(G_r)=\{I,-I\}^{r},\qquad
Dim(R_{PSL(2,\mathbb C)}(G_r))=2r.
\]
The groups $G_r$ are pairwise non-isomorphic, since their abelianizations
are $(\mathbb Z_2)^r$. Thus, for every fixed integer $t$, infinitely
many of them have zero-dimensional $SL(2,\mathbb C)$ representation
variety and $PSL(2,\mathbb C)$ representation variety of dimension
at least $t$.
In this paper we prove the following theorem quite
useful in handling many cyclically pinched one relator
presentations like those of closed non-orientable surface groups,
torus knot groups, or many of the parafree groups introduced in
\cite {B1}.

\begin{namedtheorem}{Theorem C}
 \begin{enumerate}
\item Let $G=<X,y \, ; \, W(X)=y^p>$, where
$Card(X)=n$,\, $p>1$, and $W(X)$ is a non-trivial freely reduced
word in the free group on $X$. If $A = PSL(2,\mathbb C)$, then
$Dim(R_{A}(G))=Max\{3n,Dim(R_{A}(G'))+2\}\leq 3n+1$, where $G'=<X
\, ; \, W(X)=1>$.
\item If $Dim(R_{A}(G'))+2\geq 3n$ in (1), then $R_{A}(G)$ is
reducible.
\end{enumerate}
\end{namedtheorem}

\medskip
In this paper it is also shown that if a group $G$ embeds into an
algebraic group $A$, then the profile function of $R_A(G)$, and
the profile function of $R_A(G/N(w))$, where $N(w)$ is the normal
closure of a non-trivial word $w$ in $G$, are always different.
Indeed, the representations killing $w$ form a proper closed
subset of $R_A(G)$, and the closed-inclusion lemma below shows that
its profile is strictly smaller in descending lexicographic order. This result particularly applies
to fg non-abelian fully residually free groups, a class that is
known to contain all but a finite number of fundamental groups of
closed orientable and non-orientable surfaces. In passing,
observe that if a group $G$ does
 not embed into an algebraic group $A$, then the
 profile function need not change when going to a non-trivial
 quotient
 of $G$; a good example is
the Baumslag--Solitar group
$G=\langle x,y\mid x^{-1}y^2x=y^3\rangle$.
This finitely generated group is non-Hopfian \cite{BS}, hence is not
residually finite, since finitely generated residually finite groups
are Hopfian. Choose a nonidentity element $w$ in the intersection of
the kernels of all homomorphisms from $G$ to finite groups.
For every homomorphism $\rho:G\rightarrow A$, its image is a finitely
generated linear group and is therefore residually finite by
Mal'cev's theorem \cite{MA}. It follows that $\rho(w)=I$: otherwise
a finite quotient of $\rho(G)$ would detect $w$, a contradiction.
Consequently every representation kills the normal closure $N(w)$,
and precomposition identifies $R_A(G/N(w))$ with $R_A(G)$ itself.
In particular, their profiles are equal.

\medskip
Given the positive tone of the foregoing, it is but prudent one
mentions that the next section begins with the rather sobering
result: {\it \lq\lq Given any algebraic group $A$, there are
uncountably many non-isomorphic finitely generated groups with
isomorphic representation varieties over $A$\rq \rq.} Indeed, no
invariant
 of representation varieties over $A$  can
discriminate these groups. To assist the reader in placing the
above in a more nuanced perspective, we recall that there are
infinitely many rank two parafree groups of deviation one with
equi-dimensional non-isomorphic $SL(2,\mathbb C)$ representation
varieties; see \cite {L4}. Indeed, this seems at odds with the
commonly received understanding that  parafree groups are very
much like free groups\footnote {This can also be seen as an
indication of the sensitivity of the invariants used.}, especially
when one is told that any $n$-generated group $G$ with
$Dim(R_{SL(2,\mathbb C)}(G))=3n$ is free. In fact, if the criterion
that the deviation be one is discarded, then given a parafree
group of rank $2$, there exists infinitely many parafree groups
also of rank 2 having representation varieties of dimension
greater than any apriori chosen integer \cite {L2}.

 \section*{One:   The Profile Function}

 \medskip
 Unless specified, the term \lq\lq algebraic group" is assumed
 to mean
 a complex affine algebraic group, and all algebraic varieties will always
be complex, even if some of the results in
 the sequel hold in far more general settings, as many do.
 Having made that clear,
 given any algebraic group $A$ with identity $I$,
 and a group $G$ generated by a finite set of generators
 $X=\{x_1,\cdots,x_n\}$,
 the set $Hom(G,A)$  brings with it the structure of an algebraic
 variety \footnote {Varieties are reduced complex algebraic sets and need not be irreducible. Dimension means Krull dimension, equal here to complex dimension.}.
 Indeed, this is
 simple
 to see since $Hom(G,A)$ is nothing  but $\bigcap P_{W_i(X)}^{-1}(I)$,
 where each
 $P_{W_i(X)}^{-1}(I)$ is the fiber over $I$ pertaining to the
 regular map given by $P_{W_i(X)}:A^n \rightarrow A$; this map
 arises when one evaluates in the algebraic group $A^n$ the  $i$-th relation of a possibly
 infinite set of relators
 $\{W_1=1,\cdots, W_i=1, \cdots \}$ obtained from
 a presentation of  $G$ on $X$. This algebraic
 variety goes by the name \lq\lq the representation variety" of $G$ over $A$,
 and we denote it by $R_A(G).$
 
 \medskip
  Fortunately, the Hilbert
 Basis Theorem guarantees that only a finite number of the
 relators of $G$ on the generating set $X$ are necessary in
 defining the algebraic variety $R_A(G)$. An immediate
 consequence of this is the following:

 \begin{namedtheorem}{Lemma 1.1} Given any algebraic group $A$ there
 exists an uncountable set $S$ of pairwise non-isomorphic
 2-generated groups with the property that any two groups in $S$
 have isomorphic representation varieties over the algebraic group $A$.
 \end{namedtheorem}
 \begin{proof}[Proof]
Let $S_2$ contain one representative of each isomorphism class of
two-generated groups. This set is uncountable by Neumann's theorem
\cite{NH}. For a presentation $G=\langle x,y\mid\mathcal R\rangle$,
Noetherianity of $A^2$ implies that the common zero set of all relators
in $\mathcal R$ is already defined by a finite subset
$\mathcal R_0\subseteq\mathcal R$. Thus $R_A(G)$ is the representation
variety of the finitely presented group
$\langle x,y\mid\mathcal R_0\rangle$.
There are only countably many finite subsets of the countable free
group $F_2$, and hence only countably many possible isomorphism types
of these varieties. Partition $S_2$ according to this type. Since a
countable union of countable sets is countable, at least one part is
uncountable, as required.
\end{proof}

\medskip Given the negativity of the above, it is fortunate that
one can easily show that an $n$-generated group $G$ is
free of rank $n$ iff $Dim(R_A(G))=nDim(A)$, where $A$ is an irreducible algebraic
group containing a free group of rank 2. This was proven in \cite
{L4} for $SL(2,\mathbb C)$, but the proof can be readily generalized
to an arbitrary irreducible affine group containing a free group
of rank 2.

\medskip In fact, one can do better; we next show that even when the
algebraic group $A$ is reducible, the representation variety of a
fg free group $F_n$ of rank $n$ over $A$ has remarkable properties
provided that $A$ contains a group isomorphic to $F_2$, the free
group of rank 2. To this end it will be necessary to introduce
what we shall call the \lq\lq Profile Function $P_m$" from the
space $\mathfrak S_m$ of all algebraic varieties of dimension at most
$m$ over $\mathbb C$ to $\mathbb Z_{\geq 0}^{m+1}$ given by:

\begin{equation*}
P_m(V)=(N_m(V),N_{m-1}(V),\cdots,N_0(V)),
\end{equation*}
 where $V$ is in $\mathfrak S_m$,
and $N_i(V)$ is the function that counts the number of irreducible
components  of dimension precisely $i$ in $V$.
As will become amply clear in the sequel, this invariant will
prove to be quite useful in the study of finitely generated
groups, but first it is important that we introduce some
additional notions surrounding it. Before doing so, let's make it
clear that the profile function for an algebraic variety of
dimension $m$ can be taken at values $n\ge m$ by retaining
$n-m$ leading zeros. Throughout a comparison, the same index is used.

\begin{nameddefinition}{Definition 1.1 (The Profile of a FG Group Over $A$)}If
$G$ is a finitely generated group, then the {\it profile of $G$}
over an affine algebraic group $A$ is defined to be $P_d(R_A(G))$,
where $d=Dim(R_A(G)).$
\end{nameddefinition}

\begin{nameddefinition}{Definition 1.2 (Same Profile)} Two algebraic varieties $V$
and $W$ have the {\it same profile} if they have the same
dimension  $d$, and $P_d(V)=P_d(W)$.
\end{nameddefinition}

\medskip\par\noindent

A trivial consequence of Lemma 1.1 is that given any algebraic
group $A$, there exists an uncountable set $S_2$ of pairwise
non-isomorphic 2-generated groups having the property that for
each pair of different groups $G_i,\,\,G_j \in S_2$ the profile of
the representation varieties of $G_i$ and $G_j$ over $A$,
respectively, are the same.

\begin{nameddefinition}{Definition 1.3 (Proper Relative Descent, Relative Descent)}
Let $V$ and $W$ be algebraic varieties and let $d=Dim(V)$.
We say that $W$ descends relative to $V$ if $Dim(W)\leq d$ and
there is a regular map $i:W\rightarrow V$ such that

\begin{equation*}
P_d(\overline{i(W)})=P_d(W),
\end{equation*}

where the closure is taken in $V$. The descent is proper if, in addition,
$\overline{i(W)}\subsetneq V$.
\end{nameddefinition}

This formulation excludes dense proper images. For example, the inclusion
$\mathbb C^\times\rightarrow\mathbb A^1$ preserves the profile and has proper
image, but its image is dense, so it does not give proper descent.
For a quotient $G/N$ the natural map $R_A(G/N)\rightarrow R_A(G)$ is
a closed embedding: its image consists of the representations that kill
every element of $N$. Thus it gives relative descent.

\begin{namedtheorem}{Lemma 1.0 (Closed-inclusion lemma)}
Let $Y\subsetneq X$ be a proper closed subvariety and let $d=Dim(X)$.
Then $P_d(Y)<_{\rm lex}P_d(X)$, where the coordinates are compared
from dimension $d$ down to dimension zero. The empty variety is assigned
the zero profile.
\end{namedtheorem}
\begin{proof}[Proof]
Let $k$ be the largest dimension of an irreducible component of $X$
not contained in $Y$. All components of larger dimension are retained
as components of $Y$. A missing component of dimension $k$ meets $Y$
in a proper closed subset of dimension less than $k$. A new component
of $Y$ cannot lie inside a retained component of $X$, since it would
not be maximal in $Y$. It follows that $N_i(Y)=N_i(X)$ for $i>k$ and
$N_k(Y)<N_k(X)$, proving the assertion.
\end{proof}

In particular, proper relative descent changes the profile. The inequality
is lexicographic, not coordinatewise: a finite set inside an irreducible
surface has more zero-dimensional components than the surface.

\begin{nameddefinition}{Definition 1.4 (Proper Relative Falling, Relative Falling)}
An algebraic variety $T$ is said to have {\it properly falling
profile relative to} an algebraic variety $V$ if there exists an
algebraic variety $W$ descending properly relative to $V$ with
$c=Dim(T)=Dim(W)$ and such that $P_c(T)=P_c(W)$.
\end{nameddefinition}

\medskip\par\noindent

$T$ is said to be {\it falling relative to} $V$ if the algebraic
variety $W$ above only descends relative to $V$,
$c=Dim(T)=Dim(W)$, and $P_c(T)=P_c(W)$.

\begin{namedtheorem}{Theorem 1.1} If the free group of rank 2 embeds into
the not necessarily irreducible algebraic group $A$, then an
$n$-generated group $G$ is free of rank $n$ provided that
$Dim(R_A(G))=c$ and $P_c(R_A(G))=P_c(R_A(F_n))$, where
$c=nDim(A).$
\end{namedtheorem}

\begin{proof}[Proof]
Choose an epimorphism $F_n\rightarrow G$ from the stated generating set.
If its kernel contains a nonidentity word $w$, choose an embedding
$F_n\rightarrow A$; such an embedding exists because $A$ contains $F_2$
(and hence free subgroups of every finite positive rank).
The corresponding point of $A^n=R_A(F_n)$ does not satisfy $w=1$.
Consequently $R_A(G)$ is a proper closed subvariety of $A^n$.
Lemma 1.0 contradicts the assumed equality of profiles. The kernel is
therefore trivial and $G\cong F_n$.
\end{proof}

\begin{namedtheorem}{Theorem 1.2}
Let the finitely generated group $G$ have a faithful representation in
the algebraic group $A$, and let $N$ be a nontrivial normal subgroup of
$G$. Then $R_A(G/N)$ descends properly relative to $R_A(G)$. In particular,
if $d=Dim(R_A(G))$, then

\begin{equation*}
P_d(R_A(G/N))<_{\rm lex}P_d(R_A(G)).
\end{equation*}

\end{namedtheorem}
\begin{proof}[Proof]
Precomposition with $G\rightarrow G/N$ identifies $R_A(G/N)$ with the
closed subvariety of representations killing $N$. Choose $w\in N$
with $w\neq 1$ and a faithful representation $\rho:G\rightarrow A$.
Then $\rho(w)\neq I$, so $\rho$ does not belong to that subvariety.
The inclusion is proper, and Lemma 1.0 gives the profile inequality.
\end{proof}

\begin{namedtheorem}{Corollary 1.1}
Let

\begin{equation*}
G_0\longrightarrow G_1\longrightarrow G_2\longrightarrow\cdots
\end{equation*}

be a sequence of epimorphisms of finitely generated groups, each of which
embeds in the fixed algebraic group $A$. Then the epimorphisms are
isomorphisms from some index onward. In particular, no infinite sequence
of epimorphisms with nontrivial kernels exists under these hypotheses.
\end{namedtheorem}
\begin{proof}[Proof]
Precomposition with the successive epimorphisms gives a descending chain
of closed subvarieties of $R_A(G_0)$. This chain stabilizes because its
Zariski topology is Noetherian. If a subsequent epimorphism had a
nontrivial kernel, Theorem 1.2 would make the corresponding inclusion
strict, a contradiction.
\end{proof}

\medskip Next a lemma is introduced which is quite useful in giving a
 lower bound for the dimension of representation varieties of
finitely generated groups.

\medskip
 \begin{namedtheorem}{Lemma 1.2} Let $R_A(G)$ be the representation variety
 of a group $G$ with a finite presentation of deficiency $d$ in an
 algebraic group $A$ of dimension $t$. Then all irreducible components
 of $R_A(G)$ have dimension greater than or equal
 to $dt$. In particular, $Dim (R_A(G)) \ge dt$.
 \end{namedtheorem}

\begin{proof}[Proof]
Write the presentation with $n$ generators and $m$ relators, so $d=n-m$.
The relators define a morphism $\phi:A^n\rightarrow A^m$ with
$R_A(G)=\phi^{-1}(I,\ldots,I)$. Complex algebraic groups are smooth;
all components of $A^n$ have dimension $nt$. Near the identity of $A^m$,
its maximal ideal is generated by $mt$ local parameters. Pulling these
back locally defines the fiber by at most $mt$ equations in the smooth
ambient variety $A^n$. The height theorem therefore bounds the
codimension of every component of the fiber by $mt$. Every such component
has dimension at least $nt-mt=dt$.
\end{proof}

\begin{namedtheorem}{Lemma 1.3} Let $A$ be an irreducible algebraic
 group of dimension $d$ and let $G$ have a presentation of
 deficiency $n-m$. If $Dim(R_A(G))=d(n-m)$, then for some
 positive integer $i,$
 
\begin{equation*}
P_{d(n-m)}(R_A(G))=(i,0,\dots,0).
\end{equation*}

 \end{namedtheorem}
 \begin{proof}[Proof] Every irreducible component has dimension $\ge d(n-m)$ by Lemma 1.2.
 But, $Dim(R_A(G))=d(n-m)$.  It follows then that there are only
 zeros in the rest of the entries of the profile function of its
 representation variety.
 \end{proof}

 \begin{namedremark}{Example} Let $G=<x,y,z \, ; \, x^py^qz^t=1>,$ where $p,q,$ and $t$ are
 integers greater than one with at least one of them greater than
 two. If $A=SL(2,\mathbb C),$ we have that $P_6(R_A(G))=(i,0,0,0,0,0,0),$ where
 $i$ is an integer greater
 than one. The dimension and reducibility assertions follow from \cite{L4};
 equidimensionality then follows from Lemma 1.2.
 It turns out, by a result in \cite {BC}, that if all of the integers
 are exactly 2 then $i=1$. As we shall see in the sequel, the
 situation changes if $A=PSL(2,\mathbb C)$.
 \end{namedremark}

 \begin{namedtheorem}{Theorem 1.3} Let $G$ be the fundamental group of either a
 closed orientable surface of genus $g\geq 1$ or a closed non-orientable surface
 of genus $g \geq 4,$ and let $A=SL(2,\mathbb C)$. Suppose that
 $N$ is a non-trivial normal subgroup in $G$. Let $Dim(R_A(G))=d$, then
 \begin{enumerate}
 \item $P_d(R_A(G))$ has a one in its left most entry and zeros in all
 other entries,
 \item $N_d(R_A(G/N))=0,$ and
 \item $P_d(R_A(G/N))\neq P_d(R_A(G))$.
\end{enumerate}
\end{namedtheorem}

\begin{proof}[Proof] By \cite{RBC}, closed orientable surface groups
of genus $g\geq1$ have irreducible representation varieties in
$SL(2,\mathbb C)$. The non-orientable result for $g\geq3$ is due
to \cite{BC}. This proves
(1). Now the orientable surface groups of genus $g=1$ embed into
$SL(2,\mathbb C)$, and by \cite {BGSP} all of the orientable surface
groups of genus $g \geq 2$  as well as all the non-orientable
surface groups of genus $g \geq 4$  embed into $SL(2,\mathbb C)$
since they are fully residually free and $SL(2, \mathbb C)$ is a
connected complex semisimple Lie group; so Theorem 1.2 yields (2)
together with (1). Part (3) is a result of (2).
\end{proof}

\begin{namedtheorem}{Theorem 1.4} 
Let $A$ be a nontrivial connected complex semisimple algebraic group, $G$ a
non-abelian fg fully residually free group, and $N$ a non-trivial
normal subgroup of $G$. Then $ P_d(R_A(G))\neq P_d(R_A(G/N))$.
\end{namedtheorem}

\begin{proof}[Proof] Theorem 1.5 of \cite{BGSP} gives a dense embedding
of $G$ into $A$. The result follows from Theorem 1.2.
\end{proof}

\begin{namedtheorem}{Proposition 1.1} Let $G$ be a finitely generated free
abelian group, $A=SL(2,\mathbb C)$, and suppose that $Dim(R_A(G))=d$,
and that $N$ is non-trivial normal subgroup of $G$. Then
$P_d(R_A(G))\neq P_d(R_A(G/N))$.
\end{namedtheorem}

\begin{proof}[Proof] Choose a free abelian basis
$\{x_1,\dots,x_n\}$ of $G$. Let $L$ be the list of the first
$n$ odd primes denoted thus: $p_1,p_2,\dots,p_n$. With each prime
$p_k$ in the list $L$ associate a two by two matrix $m_k$ with the
property that $m_k$ has entries $a_{i,j}=0$ if $i\neq j$,
$a_{1,1}=p_k$, and $a_{2,2}=\frac {1}{p_k}$. The matrix $m_k$ lies
in $A$, and it is straightforward to see that the homomorphism
$\phi : G\rightarrow A$ given by $\phi(x_i)=m_i$ is an embedding
of $G$ into $A$. The conclusion of the proposition follows from
Theorem 1.2.
\end{proof}

\section*{Two: Projection of Solutions and Word Equations}

\medskip In this section we will assume that all algebraic groups are
irreducible. If $A$ is an algebraic group over an algebraically
closed field of characteristic zero, $N$ is a closed normal
subgroup of $A$, and $\pi: A \rightarrow A/N$ is the projection
map, then it is well known that $A/N$ is an algebraic group. If
$N$ is finite, then it lies in the center of $A$ whenever $A$ is
irreducible. Note that every finite normal subgroup of an
algebraic group is necessarily closed.

\begin{nameddefinition}{Definition 2.1} Let $A$ be an algebraic group,
$N$ a finite normal subgroup of $A$, and  $\pi: A \rightarrow A/N$
the projection map. If $S\subseteq A$, then define
$\pi(S)=\{\pi(s) \, \vert \, s \in S \}$ to be the projection of
$S$ in $A/N.$
\end{nameddefinition}

\begin{namedtheorem}{Proposition 2.1} Let $S\subseteq A$ be an algebraic
set, then $\pi(S)$ is also an algebraic set.
\end{namedtheorem}

\begin{proof}[Proof] The projection map $\pi: A \rightarrow A/N$ is a
surjective open and closed regular map and thus, since $N$ is
finite, maps closed algebraic sets of $A$ to closed algebraic sets
of $A/N$.
\end{proof}

\medskip\par\noindent The next result is a standard fact in the theory
of algebraic groups. We omit the proof.

\begin{namedtheorem}{Proposition 2.2} Let $G$ and $G'$ be algebraic groups
with finite normal subgroups $H$ and $H',$ respectively. Then $G/H
\times G'/H'$ and $(G \times G')/(H \times H')$ are isomorphic as
algebraic groups.
\end{namedtheorem}

\medskip\par\noindent

Induction yields:

\begin{namedtheorem}{Corollary 2.1} Let \{$G_i$\} be a finite collection of
algebraic groups and \{$H_i$\} be corresponding finite normal
subgroups, respectively, where $i \in \{1,\dots,n\}.$ Then the
Cartesian product of the \{$G_i$\}'s modulo the Cartesian product
of the \{$H_i$\}'s is isomorphic, as an algebraic group, to the
Cartesian product of the \{$G_i/H_i$\}'s.
\end{namedtheorem}

\begin{namedremark}{Notation} Denote the canonical projection map from an
algebraic group $A$ to the algebraic group $A/H$ by $\pi$. Please
note that under the proper context $\pi$ shall also denote the
canonical projection map from $A^n$ to $(A/H)^n$. \end{namedremark}

\begin{namedtheorem}{Proposition 2.3} Let $H$ be a finite normal subgroup of
the algebraic group $A$ and let $\pi: A^n \rightarrow (A/H)^n$ be
the projection map. Then: \begin{enumerate}

\item $Dim(A/H)=Dim(A)$
\item $Dim(A^n)=Dim((A/H)^n)=nDim(A)$
\item If $V$ is an algebraic set in $A^n,$ then $Dim(V)=Dim(\pi(V))$.
\item If $V$ is an irreducible algebraic set, then $\pi(V)$ is also
irreducible.
\end{enumerate}
\end{namedtheorem}
\begin{proof}[Proof] \begin{enumerate}
\item : Under the conditions stipulated, $Dim(A/H)=Dim(A)-Dim(H)$;
but $H$ is finite, and thus has dimension zero.
\item : The dimension
of the product of a finite number of algebraic varieties is the
sum of the dimension of its factors. By (1), it is the case that
$Dim(A/H)=Dim(A)$. It follows that
$Dim(A^n)=Dim((A/H)^n)=nDim(A)$.
\item : The map $\pi$ is regular
 and has finite fiber over each point of $\pi(V)$. It follows that
$Dim(V)=Dim(\pi(V))$.
\item: If $V$ is irreducible and $\pi(V)$ is not, then the
irreducibility of $V$ is contradicted since the projection map
$\pi$ is a closed regular map.
\end{enumerate}
\end{proof}

Given an algebraic group $A$, a free group $F_n$ on
$\{x_1,\dots,x_n\}$, a non-trivial word $w=w(x_1,\dots,x_n)$ in
$F_n$, and the expression $w=h$, where $h\in A$, one can think of
$n$-tuples $(m_1,\dots,m_n)$ of elements of $A$ as solutions to
the equation $w=h$ provided that an orderly replacement in $w$ of
each $x_i$ with $m_i$  makes valid the statement
$w(m_1,\dots,m_n)=h$. This is the same as saying that there exists
a representation $\rho\in R_A(F_n)$ having the property that
$\rho(w)=h$. Both ways of looking at the situation will be used
interchangeably in the sequel. Note that it is not necessary that
all generators of $F_n$ appear in the word $w$.

 \begin{namedtheorem}{Theorem 2.1} Let $A$ be an
algebraic group and $H$ a finite normal subgroup of $A$ of order
$m$. Let $w$ be a non-trivial freely reduced word in the free
group $F_n$ with basis $\{a_1,\dots, a_n  \}$. Then the solutions
of the equation $w=1$ in $A/H$ are given by the projections to
$A/H$ of the solutions in $A$ of the equations $\bigcup_{i=1}^m
\{w=h_i \, \vert \, h_i \in H\}.$
\end{namedtheorem}

 \begin{proof}[Proof] If $(a_1,\dots,a_n)$ is a solution in $A$ to the
equation $w=h_i$, with $h_i \in H$, then it is clear that
$\pi(w(a_1,\dots,a_n))=\pi(h_i)$, where $\pi: A \rightarrow A/H$
is the projection map. Thus, $\pi(w(a_1,\dots,a_n))=1$. And so
solutions in $A$ of $w=h_i$ ($h_i \in H$) correspond to solutions
over $A/H$ of $w=1$. Now suppose that $w(b_1,\dots,b_n)=1$, where
$b_i \,\in A/H$; then each $b_i$ corresponds to a coset $c_iH$,
for some $c_i\, \in A$. In particular, $w(c_1,\dots,c_n)=t$, where
$t\in \pi ^{-1}(1)$; but $\pi ^{-1}(1)=H$.
\end{proof}

\begin{namedtheorem}{Proposition 2.4}
Let $V\subseteq(A/H)^n$ be nonempty and closed, and put $X=\pi^{-1}(V)$.
Then $V$ has falling profile relative to $X$. The falling is proper
exactly when some orbit of $H^n$ on the irreducible components of $X$
has more than one member.
\end{namedtheorem}
\begin{proof}[Proof]
The finite group $H^n$ acts on $X$ by coordinatewise translation, and
$\pi:X\rightarrow V$ is its finite quotient. We use the following
component argument, valid for any quotient of an affine variety by a
finite group $\Gamma$.

For each irreducible component $C$ of $X$, its image is closed and
irreducible, and $\pi^{-1}(\pi(C))=\bigcup_{\gamma\in\Gamma}\gamma C$.
If $\pi(C)\subseteq\pi(D)$ for another component $D$, irreducibility
implies $C\subseteq\gamma D$ for some $\gamma$, whence $C=\gamma D$.
Thus the components of the quotient correspond exactly to the
$\Gamma$-orbits of components upstairs. The finite map preserves their
dimensions.

Choose one component from each orbit and let $W$ be their union. Then
$W$ is closed in $X$ and has the same profile as $V$. Its inclusion in
$X$ gives descent, so Definition 1.4 gives falling of $V$ relative to
$X$. If an orbit has more than one component, $W\subsetneq X$, giving
proper falling. If every orbit is a singleton, the profiles of $V$ and
$X$ agree, and Lemma 1.0 rules out proper falling.
\end{proof}

\begin{namedremark}{Example}Let $A=SL(2,\mathbb C)$ and  $A/H=PSL(2,\mathbb C)$.
Let $V=R_{PSL(2,\mathbb C)}(\mathbb Z_2)$. Then it is easy to see that
$Dim(V)=2$ and $P_2(V)=(1,0,1)$.
\medskip\par\noindent
Notice, however, that $P_2(\pi^{-1}(V))=P_2(\{m \, \vert \,
m^2=\pm I,\, m \in SL(2,\mathbb C)\})$, the right side of which is
nothing but $P_2(R_{SL(2,\mathbb C)}(\mathbb Z_2)\cup\{m \, \vert \,
m\in SL(2,\mathbb C), \, tr(m)=0 \})=(1,0,2).$ Thus, $V$ has properly
falling profile relative to $\pi^{-1}(V)$.
\end{namedremark}

\section*{\texorpdfstring{Three: Representation Varieties in $PSL(2,\mathbb C)$ for a Class of One-Relator Groups}{Three: Representation Varieties in PSL(2,C) for a Class of One-Relator Groups}}

The next theorem gives a method for computing $Dim(R_A(G))$, where
$A=PSL(2,\mathbb C)$, for certain classes of one relator cyclically
pinched groups, and also gives conditions guaranteeing the
reducibility of the corresponding representation variety.

\medskip
\begin{namedtheorem}{Theorem 3.1} Put $A=PSL(2,\mathbb C)$.
\begin{enumerate}
\item Let $G=<X,y \, ; \, W(X)=y^p>$, where
$Card(X)=n$, $p>1$, and let $W(X)$ be a non-trivial freely reduced
word in the free group on $X$. Then
$Dim(R_A(G))=Max\{3n,Dim(R_A(G'))+2\}\leq 3n+1$, where $G'=<X \, ;
\,  W(X)=1>$.
\item If $Dim(R_A(G'))+2\geq 3n$ in (1), then $R_A(G)$ is
reducible.
\end{enumerate}
\end{namedtheorem}
\medskip

Before proceeding we shall adopt the following notation.

\begin{namedremark}{Notation} Given an algebraic group $\mathfrak A$, an element
$a\in \mathfrak A$, and an integer $p$, denote by $\Omega_{\mathfrak
A}(p,a)$ the closed set $\{m \, \vert \, m^p=a, m\in \mathfrak A\}$.
When the context is clear
 $\Omega_{\mathfrak A}(p,a)$ will be simply denoted as $\Omega(p,a)$.
 \end{namedremark}

\begin{namedremark}{Example}
Let $G=\langle x,y,z\mid[x,y]=z^p\rangle$, where $p\geq2$, and put
$A=PSL(2,\mathbb C)$ and $G'=\langle x,y\mid[x,y]=1\rangle$.
The commuting-pair variety $R_{SL(2,\mathbb C)}(G')$ is irreducible
of dimension four \cite{RBC}. Its image under $\pi$ retains these
properties. The other sign fiber is
\[
V=\{(M,N)\in SL(2,\mathbb C)^2:MN=-NM\}.
\]
It is a single simultaneous-conjugacy orbit, represented by
\[
D=\begin{pmatrix}i&0\\0&-i\end{pmatrix},\qquad
J=\begin{pmatrix}0&1\\-1&0\end{pmatrix}.
\]
Indeed, anticommutation forces both traces to be zero. After conjugating
$M$ to $D$, the matrix $N$ has the form
$\left(\begin{smallmatrix}0&b\\-b^{-1}&0\end{smallmatrix}\right)$
with $b\neq0$, and diagonal conjugation sends it to $J$.
The common centralizer of $D,J$ in $SL(2,\mathbb C)$ is
$\{I,-I\}$, so $V$ is irreducible of dimension three.
The two sign fibers are invariant under independent sign changes;
their finite images are closed and disjoint. By Theorem 2.1 they exhaust
$R_A(G')$, which therefore has profile $(1,1,0,0,0)$ with index four.
Consequently Theorem 3.1 gives
\[
Dim(R_A(G))=\max\{6,4+2\}=6
\]
and proves reducibility. For $p=2$, $G$ is the closed non-orientable
surface group of genus three: a torus with one boundary component
joined to a M\"obius band gives the relator $[x,y]z^2$, up to replacing
$z$ by its inverse. In this case $R_{SL(2,\mathbb C)}(G)$ is
irreducible by \cite{BC}, whereas $R_A(G)$ is reducible.
\end{namedremark}

\begin{namedtheorem}{Lemma 3.1}
For $A=PSL(2,\mathbb C)$ and $t\geq2$, $Dim(R_A(\mathbb Z_t))=2$.
\end{namedtheorem}
\begin{proof}[Proof]
The solutions of $a^t=I$ consist of the identity and finitely many
nontrivial semisimple conjugacy classes. At least one such class exists,
and each has dimension two. Theorem 4.4 gives their precise number.
\end{proof}

\begin{namedtheorem}{Lemma 3.2}
A dominant morphism from an affine variety to an irreducible variety
has an image containing a nonempty open subset.
\end{namedtheorem}
\begin{proof}[Proof] This follows from constructibility of the image; see
\cite{MD}. \end{proof}

\begin{namedtheorem}{Lemma 3.3}
Let $Z$ be a proper closed subvariety of a variety $W$. If
$Dim(W-Z)=n$ and $Dim(Z)\geq n$, then $W$ is reducible.
\end{namedtheorem}
\begin{proof}[Proof]
If $W$ were irreducible, its nonempty open subset $W-Z$ would have
dimension $Dim(W)=n$, whereas its proper closed subset $Z$ would have
dimension less than $n$.
\end{proof}

\begin{proof}[Proof of Theorem 3.1]
Throughout, $A=PSL(2,\mathbb C)$. Every element of $A$ has a $p$th root.
Indeed, a semisimple element has a root in a maximal torus. A nontrivial
unipotent element has a lift conjugate to
$\left(\begin{matrix} 1&u\\0&1\end{matrix}\right)$, with root obtained by replacing
$u$ by $u/p$. The identity also has roots. For $a\neq I$ the set
$\Omega_A(p,a)$ is finite: if $b$ is a lift of a root and $c$ a lift
of $a$, then $b^p=\pm c$ and $b$ commutes with the noncentral matrix
$c$. Its Jordan form shows that each of these equations has finitely
many solutions in the centralizer of $c$.

The projection $\Phi:R_A(G)\rightarrow A^n$ forgetting the last
coordinate is therefore surjective. Put $B=R_A(G')$ and
$Z=\Phi^{-1}(B)$. Directly from the equations,

\begin{equation*}
Z=B\times\Omega_A(p,I),\qquad Dim(Z)=Dim(B)+2.
\end{equation*}

The word $W$ is nontrivial in $F_n$, and $F_n$ embeds in $A$.
Consequently $B\subsetneq A^n$. Since $A^n$ is irreducible,

\begin{equation*}
3(n-1)\leq Dim(B)\leq3n-1;
\end{equation*}

the lower bound is Lemma 1.2.

Over the nonempty open subset $A^n-B$, all fibers of $\Phi$ are
nonempty and finite. The dimension theorem for morphisms gives

\begin{equation*}
Dim(R_A(G)-Z)=3n.
\end{equation*}

Taking the maximum of the dimensions of this open part and its closed
complement proves

\begin{equation*}
Dim(R_A(G))=\max\{3n,Dim(B)+2\}\leq3n+1.
\end{equation*}

If $Dim(B)+2\geq3n$, Lemma 3.3 proves reducibility.
\end{proof}

\section*{Four: Algebraic Varieties with Conditions}

If $A$ is an algebraic group and $t \in A,$ then the left
translation $\rho_t:A\rightarrow A$ given by $\rho_t(v)=tv$ is an
isomorphism of the algebraic variety $A$ onto itself. Next a
condition is defined ensuring that certain closed sets in $A$ are
mapped to themselves by  translations $\rho_t$ associated with any
element $t$ in the center of the algebraic group $A$. We could
define this for more general settings, but the case
$A=R_{SL(2,\mathbb C)}(F_n)$ is enough for our treatment.

\begin{nameddefinition}{Definition 4.1. (The $\pm$ Condition)}An algebraic variety
 $V\subseteq SL(2,\mathbb C)^n$ satisfies the {\it $\pm$ Condition} if whenever
 $v=(m_1,\dots,m_n)$ is in $V$, then so is $(\pm m_1,\dots,\pm
 m_n)$.
 \end{nameddefinition}

 \par\noindent So an algebraic variety $V$ satisfies the $\pm$
 Condition provided it is mapped to itself by all isomorphisms
 of  $SL(2,\mathbb C)^n$ induced by central
 translations.

\begin{namedtheorem}{Proposition 4.1}
Let $w\in F_n$ be nontrivial and let $V_1$ and $V_{-1}$ be the
fibers of its evaluation map over $I$ and $-I$ in $SL(2,\mathbb C)^n$.
Suppose both satisfy the $\pm$ Condition. Then
$\pi(V_1)\cap\pi(V_{-1})=\emptyset$. If both fibers are nonempty,

\begin{equation*}
R_{PSL(2,\mathbb C)}(\langle x_1,\ldots,x_n\mid w=1\rangle)
=\pi(V_1)\sqcup\pi(V_{-1})
\end{equation*}

is a disjoint union of two nonempty closed-and-open subsets, and is
therefore disconnected. The union $V_1\cup V_{-1}$ is also disconnected.
No assertion of disconnectedness is made about either individual fiber.
\end{namedtheorem}
\begin{proof}[Proof]
The fibers of $\pi:SL(2,\mathbb C)^n\rightarrow PSL(2,\mathbb C)^n$
are the orbits of independent sign changes. If an orbit met both
$V_1$ and $V_{-1}$, sign invariance would imply that a tuple evaluates
$w$ to both $I$ and $-I$, which is impossible. The two images are closed
because $\pi$ is finite; their union is the asserted representation
variety by Theorem 2.1. Each is open relative to that union because its
complement is the other closed subset. The same reasoning applies to
$V_1\cup V_{-1}$. Individual fibers can be connected: for $w=x^2$,
$V_{-1}$ is the irreducible trace-zero conjugacy class.
\end{proof}

\medskip
The $\pm$ Condition for an algebraic variety $V$ in $SL(2,\mathbb
C)^n$ is equivalent to the following mapping criteria stipulated
in the next lemma.

\begin{namedtheorem}{Lemma 4.1} 
An algebraic variety  $W$ in $R_{SL(2,\mathbb C)}(F_n)$ satisfies the
$\pm$ Condition iff $\pi^{-1}(\pi(W)) = W$, where $\pi : SL(2,
\mathbb C)^{n} \rightarrow PSL(2, \mathbb C)^{n}$ is the canonical
projection.
\end{namedtheorem}

\begin{proof}[Proof]
Each fiber of $\pi$ is precisely the set of all independent sign changes
of one tuple. Thus $\pi^{-1}(\pi(W))=W$ if and only if $W$ is a union
of these sign orbits, which is exactly the $\pm$ Condition.
\end{proof}

\medskip
\begin{namedremark}{Example} Consider the word $x^2$ in the free group of
rank one, and the resulting equations $x^2=I$, and $x^2=-I$ in
$SL(2,\mathbb C)$. It is not difficult to see that the variety $V_1$
obtained as solutions to the first equation is $\{\pm I\}$, and
that $V_1$ satisfies the $\pm$ Condition; $V_1$ is a reducible
variety with two points as irreducible components; the individual
components are permuted by the isomorphism of $SL(2,\mathbb C)$
induced by left translation  with the element $-I.$ Now, let
$V_{-1}$ be the solutions to the second equation $x^2=-I$; then
$V_{-1}$ consists of all matrices of trace zero, and since
$V_{-1}$ is the orbit under conjugation of a matrix of trace other
than $\pm 2$ , it is an irreducible variety of $SL(2,\mathbb C)$.
Notice also that the variety satisfies the $\pm$
Condition, and that left translation by the element $-I$
maps this variety to itself. Observe that $\pi(V_1)$ is the
identity element in $PSL(2, \mathbb C)$, and $\pi(V_{-1})$ are the
equivalency classes of all matrices of trace zero in $PSL(2,\mathbb
C)$. So no component of $V_1$ mapped into a component of $V_{-1}$,
and no component of $V_{-1}$ mapped into a component of $V_1$.
The two point components of $V_1$ have the same image, so
$\pi(V_1)$ is a single point.
\end{namedremark}

\begin{namedtheorem}{Theorem 4.1} 
Let $w$ be a non-trivial freely reduced word in the commutator
subgroup of $F_n=<x_1,\dots,x_n>$, and $G=<x_1,\dots,x_n \, ; \,
w=1>$. If $R_{SL(2,\mathbb C)}(G)$ is an irreducible variety and
$V_{-1}=\{\rho \, \vert \, \rho \in R_{SL(2,\mathbb
C)}(F_n),\rho(w)=-I\}\neq \emptyset$, then $P_d(R_{SL(2,\mathbb
C)}(G))\neq P_d(R_{PSL(2,\mathbb C)}(G)).$
\end{namedtheorem}

\begin{proof}[Proof]
The exponent sum $e_j(w)$ on each generator is zero, since $w$ belongs
to the commutator subgroup. For signs $\epsilon_j\in\{1,-1\}$,

\begin{equation*}
w(\epsilon_1m_1,\ldots,\epsilon_nm_n)
=\left(\prod_j\epsilon_j^{e_j(w)}\right)w(m_1,\ldots,m_n)
=w(m_1,\ldots,m_n).
\end{equation*}

Both fibers therefore satisfy the $\pm$ Condition. The image
$\pi(V_1)$ is irreducible because $V_1=R_{SL(2,\mathbb C)}(G)$ is
irreducible. Proposition 4.1 gives a disjoint nonempty second closed
piece $\pi(V_{-1})$. The projective representation variety is thus
reducible, whereas $V_1$ is irreducible, so the profiles differ.
\end{proof}

\begin{namedtheorem}{Corollary 4.1} 
Let $G$ be the fundamental group of a closed orientable surface
of genus $g \geq 1$, then $P_d(R_{SL(2,\mathbb C)}(G))\neq
P_d(R_{PSL(2,\mathbb C)}(G))$,\,where $d$ is the maximum of the dimensions of these two representation varieties.
\end{namedtheorem}

\begin{proof}[Proof]A presentation of such a group can be found with its
relator $w$ in the commutator subgroup of the free group on its
$2g$ generators; thus we can obtain two algebraic varieties $V_1$,
and $V_{-1}$ in
$SL(2,\mathbb C)^{2g}$ satisfying the $\pm $ Condition by the exponent-sum calculation in Theorem 4.1.  
Now, by \cite {RBC} we know that $R_{SL(2,\mathbb C)}(G)$ is
irreducible. The result follows once we show that $V_{-1}$ is not
empty. For example, take
$m_1=\operatorname{diag}(i,-i)$ and
$m_2=\left(\begin{smallmatrix}0&1\-1&0\end{smallmatrix}\right)$.
Then $[m_1,m_2]=-I$, where $[m_1,m_2]$ stands for the commutator
of the matrices $m_1$ and $m_2$. It follows that a word in the free
group $F_{2g}$ of the type
$[x_1,x_2][x_3,x_4]\dots[x_{2g-1},x_{2g}]$ can be evaluated over
$SL(2,\mathbb C)^{2g}$ to equal $-I$ simply by mapping each of the
generators $x_i$ for $i<2g-1$ to $I$ and setting $x_{2g-1}=m_1$
and $x_{2g}=m_2$. So the algebraic variety $V_{-1}$ is not empty;
the result follows.
\end{proof}

\begin{namedtheorem}{Theorem 4.2} 
Let $w$ be a non-trivial and freely reduced word in the free group
$F_n$, with even exponent sum on each generator. Let
  $G=<x_1,\dots,x_n; \, w=1>$. If $R_{SL(2,\mathbb C)}(G)$ is an
  irreducible variety and
$V_{-1}=\{\rho \, \vert \, \rho \in R_{SL(2,\mathbb
C)}(F_n),\rho(w)=-I\}\neq \emptyset$, then $P_d(R_{SL(2,\mathbb
C)}(G))\neq P_d(R_{PSL(2,\mathbb C)}(G)).$
\end{namedtheorem}

\begin{proof}[Proof]
For independent signs $\epsilon_j$, the sign factor in the identity
of Theorem 4.1 is $\prod_j\epsilon_j^{e_j(w)}=1$, since every
$e_j(w)$ is even. Hence both fibers satisfy the $\pm$ Condition.
Proposition 4.1 decomposes the projective representation variety
into the disjoint nonempty closed subsets $\pi(V_1)$ and
$\pi(V_{-1})$. It is reducible, while $V_1$ is irreducible by
hypothesis, so their profiles differ.
\end{proof}

The next lemma guarantees that $V_{-1}$ is always non-empty
provided that a minor assumption on the word $w$ is made.

\begin{namedtheorem}{Lemma 4.2} Let a freely reduced word $w$ in $F_n=<x_1,\dots,x_n>$
be such that it has even exponent sum not
equal to zero on one of the generators appearing in it, then the
algebraic variety $V_{-1}=\{\rho \, \vert \, \rho \in R_{SL(2,\mathbb
C)}(F_n),\, \rho(w)=-I \}\neq \emptyset$. If, in addition,
$V_{-1}$ and $V_{1}$ each satisfy the $\pm$ Condition, then $Dim
(V_{-1})\geq 3(n-1)$.
\end{namedtheorem}

\begin{proof}[Proof]
Let the nonzero even exponent sum on $x_i$ be $s$. Choose
$\lambda\in\mathbb C^\times$ with $\lambda^s=-1$ and set
$m=\operatorname{diag}(\lambda,\lambda^{-1})$. Sending $x_i$ to $m$
and all other generators to $I$ gives a point of $V_{-1}$.
If the two fibers satisfy the $\pm$ Condition, Proposition 4.1 shows
that $\pi(V_{-1})$ is a nonempty closed-and-open union of components
of $R_{PSL(2,\mathbb C)}(G)$. Lemma 1.2 bounds each of those components
below by $3(n-1)$. Since $\pi$ is finite,
$Dim(V_{-1})=Dim(\pi(V_{-1}))\geq3(n-1)$.
\end{proof}

\begin{namedtheorem}{Corollary 4.2}  
Let $G$ be the fundamental group of a closed non-orientable
surface of genus $g \geq 3$, then $P_d(R_{SL(2,\mathbb C)}(G))\neq
P_d(R_{PSL(2,\mathbb C)}(G))$,\,where $d$ is the maximum of the dimensions of these two representation varieties.
\end{namedtheorem}

\begin{proof}[Proof]
Use the closed non-orientable surface presentation with relator
$w=x_1^2\cdots x_g^2$. Both sign fibers are invariant under independent
sign changes. Lemma 4.2 shows that $V_{-1}$ is nonempty. By \cite{BC},
$V_1=R_{SL(2,\mathbb C)}(G)$ is irreducible for $g\geq3$.
Theorem 4.2 now applies.
\end{proof}

\begin{nameddefinition}{Definition 4.2 (Lift of a Representation)}
Let $\rho \in R_{PSL(2, \mathbb C)}(G)$ be a representation of a
finitely generated group $G$. Then $\rho$ is said to have a {\it
lift} if there exists $\phi \in R_{SL(2, \mathbb C)}(G)$ such that
$\pi \circ \phi=\rho$.
\end{nameddefinition}

\begin{namedremark}{Observation}  
If a group $G$ is generated by the set $\{x_1,\dots,x_n\}$, then
$R_A(G)$ is in one-to-one correspondence with the set $
\{(\rho(x_1),\dots,\rho(x_n)) \in A^n \, \vert \, \rho \in
R_A(G)\}$.
\end{namedremark}

\begin{namedtheorem}{Theorem 4.3}  
Let $G$ be a one relator $n$-generated group whose relator is a
non-trivial word $w$ in $F_n=<x_1,\dots,x_n>$ which gives rise to
algebraic varieties $V_1$ and $V_{-1}$, as in Proposition 4.1,
satisfying the $\pm$ Condition. Then no representation $\rho \in
R_{PSL(2,\mathbb C)}(G)$ with $(\rho(x_1),\dots,\rho(x_n))\in
\pi(V_{-1})$ has a lift.
\end{namedtheorem}

\begin{proof}[Proof]
If a projective representation lies in $\pi(V_{-1})$, every tuple
lifting its generator images belongs to $V_{-1}$ by sign invariance.
A homomorphism $G\rightarrow SL(2,\mathbb C)$ would instead give a tuple
in $V_1$, since it must satisfy $w=I$. This is impossible.
\end{proof}

\begin{namedtheorem}{Lemma 4.3}  
Let $V_1$ and $V_2$ be nonempty algebraic varieties in
$SL(2,\mathbb C)^n$ which satisfy the \lq\lq Minus Condition"
stipulated below.

\par\noindent The Minus Condition:

\begin{enumerate}
\item If $m=(m_1,m_2,\dots,m_n)\in V_1 $ then $-m \in V_2$.
\item If $m=(m_1,m_2,\dots,m_n)\in V_2 $ then $-m \in V_1$.
\end{enumerate}

\par\noindent Then $V_1$ and $V_2$ are isomorphic,
$\pi(V_1)=\pi(V_2)$, and if either $V_1$ or $V_2$ is irreducible,
then $\pi(V_1 \cup V_2)$ is irreducible in $PSL(2,\mathbb C)^n$.
\end{namedtheorem}

\begin{proof}[Proof] Clearly, two algebraic varieties meeting the Minus
Condition are isomorphic since the condition gives rise to a
bijective correspondence originating from a regular map whose
inverse is also a regular map. That $\pi (V_1)=\pi(V_2)$ follows
from the fact that for each point $m \in V_1$, the point $-m$ is
in $V_2$, and conversely, and the pair of points $m$, $-m$ map to
a single point in $PSL(2,\mathbb C)^n$. The irreducibility follows
since the image of an irreducible variety under a regular map
should also be irreducible.
\end{proof}

\begin{namedremark}{Example} Let $G=<x,y,z \, ; \, x^2=1>$. It is easy to see
that $R_{SL(2,\mathbb C)}(G)=\{I\times SL(2,\mathbb C)^2\}\cup
\{-I\times SL(2,\mathbb C)^2\}$. Let $V_1=I\times SL(2,\mathbb C)^2$,
and let $V_2= -I\times SL(2,\mathbb C)^2$. Then $V_1$ and $V_2$ meet
the Minus Condition, and at least one of the $V_i$ ($i = 1,2$) is
irreducible.  Consequently, $\pi(R_{SL(2,\mathbb C)}(G))$ is
irreducible. In fact, since $R_{SL(2,\mathbb C)}(G)$ is six
dimensional, with  irreducible components $V_1$, and $V_2$
respectively, $P_6(R_{SL(2,\mathbb C)}(G))=(2,0,0,0,0,0,0)$. Further,
as a consequence of the earlier statements in this example,
$P_6(\pi(R_{SL(2,\mathbb C)}(G)))=(1,0,0,0,0,0,0)$. It is readily
seen, however, that $R_{PSL(2,\mathbb C)}(G)$ consists of more than
$\pi(R_{SL(2,\mathbb C)}(G))$. In fact it contains an eight
dimensional irreducible variety, namely $\pi(V_3)$, where
$V_3=\{(m_1,m_2,m_3)\, \vert \, m_i\in SL(2,\mathbb C), tr(m_1)=0\}$.
Now, by Proposition 4.1, it is the case that since the word $x^2$
is a non-trivial word in the free group of rank 3 on the
generators $\{x,y,z\}$, and that  since  the $SL(2,\mathbb C)$
solutions $S_1$, and $S_{-1}$ to the equations $x^2=I$, and
$x^2=-I$, respectively, satisfy the $\pm$ Condition of Proposition
4.1, that $\pi(S_1 \cup S_{-1})$ is a reducible variety. So
$\pi(S_i)$ does not collapse onto $\pi(S_j)$, where $i\neq j$ and
$i,j \in \{1,-1\}$. In fact, $\pi(S_{-1})$ is irreducible since
$S_{-1}$ is irreducible and the projection map is regular. Putting
the foregoing together yields for example that $P_8(R_{PSL(2,\mathbb
C)}(G))=(1,0,1,0,0,0,0,0,0)$.
\end{namedremark}

\medskip
Next we compute the profile function of $R_{PSL(2,\mathbb C)}(\mathbb
Z_n)$, when $n\ge 2$. Clearly, by Lemma 3.1, $Dim(R_{PSL(2,\mathbb
C)}(\mathbb Z_n))=2$. Before starting it is prudent to introduce some
notational convention.

\begin{namedremark}{Notation} \begin{enumerate}

\item If $x \in SL(2,\mathbb C)$, denote by $[x]$ the equivalency class
corresponding to $x$ in $PSL(2,\mathbb C)$.
\item For integers $p,q$ and a fixed positive integer $n$, let $diag(p,q)$ stand for
the $2 \times 2$ matrix \smallskip $ \left(\begin{matrix}
e^{\frac {p \pi i} {n} } & 0\\
0 & e^{\frac {q \pi i} {n}}
\end{matrix}\right).$ Exponents may be read modulo $2n$; in the
 determinant-one matrices below we use $q=-p$. \smallskip
\item $x\sim y$ shall mean that $x$ is conjugate to $y$ (i.e. $x$
and $y$ are in the same orbit under the action of conjugation).
\end{enumerate}
\end{namedremark}

\begin{namedtheorem}{Theorem 4.4} Let $n\geq2$.
\begin{enumerate}
\item $P_2(R_{PSL(2,\mathbb C)}(\mathbb Z_n))=(\frac{n} {2},0,1)$, if $n$ is even.
\item $P_2(R_{PSL(2,\mathbb C)}(\mathbb Z_n))=(\frac{n-1} {2},0,1)$, if $n$ is odd.
\end{enumerate}
\end{namedtheorem}
\begin{proof}[Proof]
Write $\zeta=e^{\pi i/n}$ and
$D_p=\operatorname{diag}(\zeta^p,\zeta^{-p})$, with indices modulo
$2n$. A lift $M$ of an element of $R_{PSL(2,\mathbb C)}(\mathbb Z_n)$
satisfies $M^n=I$ or $M^n=-I$. Its minimal polynomial divides the
separable polynomial $T^{2n}-1$, so it is diagonalizable and conjugate
to some $D_p$.

The indices $p=0,n$ give $I,-I$ and one point downstairs. For all
other indices, the conjugacy class of $D_p$ is the closed trace fiber
of trace $\zeta^p+\zeta^{-p}\neq\pm2$. It is irreducible of dimension
two: its centralizer is the one-dimensional diagonal torus.
The conjugacy classes of $D_p$ and $D_q$ have the same projective
image precisely when
\[
q\equiv p,\ -p,\ p+n,\text{ or }-p+n\pmod{2n}.
\]
Indeed, projective conjugacy lifts to conjugacy with either $D_q$ or
$-D_q$. First identify $p$ with $-p$, leaving the representatives
$1,\ldots,n-1$. Negation of matrices then identifies $p$ with $n-p$.
This involution has no fixed point when $n$ is odd and one fixed
point $p=n/2$ when $n$ is even. It therefore has $\lfloor n/2\rfloor$
orbits.

Their images under the finite map $\pi$ are distinct closed
irreducible varieties of dimension two. Together with the identity
point they are a finite disjoint closed cover, hence exactly the
irreducible components. There are no components of dimension one.
Thus the profile is $(\lfloor n/2\rfloor,0,1)$, as asserted.
\end{proof}

Compare Theorem 4.4 
with the following:

\begin{namedtheorem}{Proposition 4.2} Let $n\geq1$.
\begin{enumerate}
\item $P_2(R_{SL(2,\mathbb C)}(\mathbb Z_n))=(\frac{n-2} {2},0,2)$, if $n$ is even.
\item $P_2(R_{SL(2,\mathbb C)}(\mathbb Z_n))=(\frac{n-1} {2},0,1)$, if $n$ is odd.
\end{enumerate}
\end{namedtheorem}

\begin{proof}[Proof]
Every solution of $M^n=I$ is diagonalizable. Its eigenvalues are
$\lambda,\lambda^{-1}$ with $\lambda^n=1$, and inversion identifies
exactly the pairs defining the same conjugacy class. Each noncentral
class is a closed irreducible surface, as in the proof of Theorem 4.4.
The central solutions are $I,-I$ for even $n$ and just $I$ for odd
$n$. Removing them leaves $n-2$ or $n-1$ roots, respectively, paired
by inversion with no fixed points. This gives the displayed counts.
\end{proof}

\medskip

We will employ  Theorem 4.4  
in computing the profile function for $R_{PSL(2,\mathbb C)}(G)$,
where $G$ is the free product of two finite cyclic groups; this
will prove useful in our study of the profile functions  of
$R_{PSL(2,\mathbb C)}(G)$, when $G$ is a torus knot group. The
next lemma will prove indispensable.

\begin{samepage}
\begin{namedtheorem}{Lemma 4.4}  
 Let $V$ and $W$ be nonempty algebraic varieties with $dim(V)=d_1$ and $dim(W)=d_2.$
 Then
\begin{enumerate}
\item $dim(V \times W)=d_1+d_2$ and
\item $N_{d_1+d_2}(V \times W)=N_{d_1}(V)N_{d_2}(W)$.
\end{enumerate}
\end{namedtheorem}
\end{samepage}

\begin{proof}[Proof] This is an elementary result, the proof of which may
be deduced from the basic properties of irreducible varieties and
their Krull dimension. See, for example, \cite {MD}.
\end{proof}

\medskip

The next lemma is stated merely for completion, and its proof will
also be omitted, as it is an elementary exercise. As usual,
$G_1*G_2$ denotes the free product of the groups $G_1$ and $G_2.$

 \begin{namedtheorem}{Lemma 4.5}  
If $G_1$ and $G_2$ are fg groups and $A$ is an algebraic group,
then $R_A(G_1*G_2)\cong R_A(G_1) \times R_A(G_2)$ as algebraic varieties.
\end{namedtheorem}

The next theorem  gives  the profile function for the $PSL(2,\mathbb
C)$ representation variety of the free product of two cyclic
groups of finite order. Besides being an interesting result in its
own right, it will also be indispensable in our study of the
profile function over $PSL(2,\mathbb C)$ for the representation
varieties associated with a class of one relator groups containing
the class of torus knot groups.

\begin{namedtheorem}{Theorem 4.5}  
Let $m,n\geq2$, let $\mathbb Z_k$ denote the cyclic group of
order $k$, and put $A=PSL(2,\mathbb C)$. Then $P_{4}(R_A(\mathbb Z_{m}*\mathbb Z_{n}))$
equals: \begin{enumerate}
\item $(\frac{mn}{4}, 0, \frac{m+n}{2}, 0, 1)$
if both $m$ and $n$ are even.

\item $(\frac{n(m-1)}{4}, 0, \frac{m+n-1}{2}, 0, 1)$
if $m$ is odd and $n$ is even.

\item $(\frac{m(n-1)}{4}, 0, \frac{m+n-1}{2}, 0, 1)$
if $m$ is even and $n$ is odd.

\item $(\frac{(m-1)(n-1)}{4}, 0, \frac{m+n-2}{2}, 0, 1)$
if both $m$ and $n$ are odd.
\end{enumerate}

\end{namedtheorem}

\begin{proof}[Proof]
Put $a=\lfloor m/2\rfloor$ and $b=\lfloor n/2\rfloor$.
By Theorem 4.4 the two cyclic representation varieties have, respectively,
$a$ and $b$ two-dimensional components, one zero-dimensional component
each, and no other components. Lemma 4.5 identifies the free-product
representation variety with their product. Its irreducible components
are exactly products of components of the factors. Thus

\begin{equation*}
N_4=ab,\qquad N_2=a+b,\qquad N_0=1,\qquad N_3=N_1=0.
\end{equation*}

Substituting the two possible parities of $m$ and $n$ gives the four
displayed profiles.
\end{proof}

Next a theorem is proven, a special case of which gives the number
of four dimensional  irreducible components of a torus knot
group. The theorem is a consequence of Theorem 3.1, Theorem 4.4,   
Lemma 4.4, 
and some of the ideas developed in \cite {L3} once the necessary
adjustments for the projective linear group are made. An immediate
consequence, when both exponents are at least three and relatively prime,
is that $N_4(R_{PSL(2,\mathbb C)}(G))$ is strictly smaller than the genus
of the torus knot, whereas $N_4(R_{SL(2,\mathbb C)}(G))$ equals that genus.

\begin{namedtheorem}{Proposition 4.3} Let $p,\,t$ be integers greater than one,
$G_{pt} = <x,y \, | \, x^{p}=y^{t}>$, and $A=PSL(2,\mathbb C)$. Then
$N_{4}(R_A(G_{pt}))$ equals:

\begin{enumerate}
\item $\frac{(p-1)(t-1)}{4}$ if both $p$ and $t$ are odd.
\item $\frac{p(t-1)}{4}$ if $p$ is even and $t$ is odd.
\item $\frac{t(p-1)}{4}$ if $p$ is odd and $t$ is even.
\item $\frac{pt}{4}$ if both $p$ and $t$ are even.
\end{enumerate}
\end{namedtheorem}

\begin{proof}[Proof]
Theorem 3.1 and Lemma 3.1 give $Dim(R_A(G_{pt}))=4$.
Inside this variety put
\[
S=\{(M,N)\in R_A(G_{pt}):M^p=I\}
 =\Omega_A(p,I)\times\Omega_A(t,I).
\]
Projection to $M$ maps $R_A(G_{pt})-S$ onto the nonempty open subset
$A-\Omega_A(p,I)$. Its fibers are nonempty and finite, by the
power-map argument in the proof of Theorem 3.1. Hence
$Dim(R_A(G_{pt})-S)=3$.
Every four-dimensional component of $R_A(G_{pt})$ must therefore lie
in $S$. Conversely, every four-dimensional component of the closed
subset $S$ is a component of the four-dimensional ambient variety.
By Theorem 4.4 and Lemma 4.4,
\[
N_4(R_A(G_{pt}))=N_4(S)
=\lfloor p/2\rfloor\lfloor t/2\rfloor.
\]
Expanding according to parity gives all four formulas.
\end{proof}

\medskip

For the sake of comparison, we state the result established in
\cite {L3} for the groups $G_{pt}$ as in the above theorem.

\begin{namedtheorem}{Theorem 4.6 (\cite {L3})} Let $G_{pt}=<x,y \, ; \,
x^p=y^t>$, where $p,t$ are integers greater than one, and
$A=SL(2,\mathbb C)$. Then \begin{enumerate} \item\, $N_4( R_A(G_{pt}))=\frac
{(p-2)(t-2)+pt} {4}$ if both $p$, $t$ are even and
\item\, $N_4( R_A(G_{pt}))=\frac { (p-1)(t-1)} {2}$ if either
$p$ or $t$ is odd.
\end{enumerate}
\end{namedtheorem}

The next theorem follows from Proposition 4.3 and Theorem 4.6.

\begin{namedtheorem}{Theorem 4.7}  
If $p,t \geq 3$ are integers, then $P_4(R_{SL(2,\mathbb
C)}(G_{pt}))\neq P_4(R_{PSL(2,\mathbb C)}(G_{pt})).$
\end{namedtheorem}

\begin{proof}[Proof]
Subtract the $PSL(2,\mathbb C)$ count in Proposition 4.3 from the
$SL(2,\mathbb C)$ count in Theorem 4.6. The difference is

\begin{equation*}
\begin{cases}
(p-1)(t-1)/4,&p,t\text{ odd},\\
(p-2)(t-1)/4,&p\text{ even},\ t\text{ odd},\\
(p-1)(t-2)/4,&p\text{ odd},\ t\text{ even},\\
(p-2)(t-2)/4,&p,t\text{ even}.
\end{cases}
\end{equation*}

Every expression is positive for $p,t\geq3$, so the dimension-four
entries, and hence the profiles, differ.
\end{proof}

Every torus knot group has a presentation of the type $<x,y \, ;
\, x^p=y^t>$, where $p$ and $t$ are relatively prime integers.
Immediate consequences of Theorem 4.7  
are the following:

\begin{namedtheorem}{Corollary 4.3}  
If a torus knot group has presentation $<x,y \, ; \, x^p=y^t>$,
where $p,t \geq 3$ are relatively prime integers, then
$N_4(R_{PSL(2,\mathbb C)}(G_{pt}))$ is not equal to the genus of the
torus knot.
\end{namedtheorem}

\begin{proof}[Proof] $N_4(R_{SL(2,\mathbb C)}(G_{pt}))$, in such a case of
the $p$ and $t$, is equal to the genus of the torus knot, \cite
{L3}.
By Theorem 4.7,  
$N_4(R_{SL(2,\mathbb C)}(G_{pt}))\neq N_4(R_{PSL(2,\mathbb
C)}(G_{pt}))$.
\end{proof}

\begin{namedtheorem}{Corollary 4.4}  
Let $G$ be a torus knot group with presentation $<x,y \, ; \,
x^p=y^t>$ with $p,t \geq 3$, then $P_4(R_{SL(2,\mathbb C)}(G))\neq
P_4(R_{PSL(2,\mathbb C)}(G))$.
\end{namedtheorem}

\begin{proof}[Proof] This is a direct consequence of Theorem 4.7 
since, for such values of $p$ and $t,$ the number of four
dimensional irreducible components are not the same. So the profile function does
not agree at dimension four.
\end{proof}

\bigskip
\begin{namedtheorem}{Lemma 4.6}  
Let $G$  be a one relator n-generated group with non-trivial
relator $w=1$. Let $V_{-1}$ and  $V_{1}$ be all representations of the free group
$F_n$ in $SL(2,\mathbb C)$ mapping $w$ to $-I$, and $w$ to $I$,
respectively. If $V_1$ and $V_{-1}$ satisfy the Minus Condition,
then every representation $\rho \in R_{PSL(2,\mathbb C)}(G)$ has a
lift.
\end{namedtheorem}

\begin{proof}[Proof]
Choose lifts in $SL(2,\mathbb C)$ of the generator images of the given
projective representation. Their tuple belongs to $V_1\cup V_{-1}$.
If it lies in $V_{-1}$, simultaneous negation moves it into $V_1$ by
the Minus Condition without changing its projective image. In either
case we obtain a tuple satisfying $w=I$, hence a lifting homomorphism.
\end{proof}

We finish with a comparison result.

\begin{namedtheorem}{Theorem 4.8}  
 Suppose that $w$ in $F_n=<x_1,\dots,x_n>$
 is a freely reduced and non-trivial word, and $G=<x_1,\dots,x_n \, ; \, w=1>$.
 Suppose that $V_{-1}=\{\rho \, \vert \, \rho \in R_{SL(2,\mathbb C)}(F_n), \rho(w)=-I\}$, and
 $V_{1}=\{\rho \, \vert \, \rho \in R_{SL(2,\mathbb C)}(F_n), \rho(w)=I\}$ are
 algebraic varieties that satisfy the Minus Condition. Then
 $R_{PSL(2,\mathbb C)}(G)$ falls, or falls properly relative to
 $R_{SL(2,\mathbb C)}(G)$.
 \end{namedtheorem}

\begin{proof}[Proof]
Let $e_j$ be the exponent sum of $x_j$ in $w$, and let

\begin{equation*}
\Gamma_w=\{(\epsilon_1,\ldots,\epsilon_n)\in\{1,-1\}^n:
\prod_j\epsilon_j^{e_j}=1\}.
\end{equation*}

The identity

\begin{equation*}
w(\epsilon_1m_1,\ldots,\epsilon_nm_n)
=\left(\prod_j\epsilon_j^{e_j}\right)w(m_1,\ldots,m_n)
\end{equation*}

shows that $\Gamma_w$ preserves $V_1$. By Lemma 4.6 the finite map
$\pi:V_1\rightarrow R_{PSL(2,\mathbb C)}(G)$ is surjective, and its
fibers are precisely the $\Gamma_w$-orbits. The component-orbit argument
in Proposition 2.4 therefore applies: select one irreducible component
from each orbit. Their union is a closed subvariety of $V_1$ having
the profile of the projective representation variety. It is proper
exactly when some orbit contains more than one component. This proves
the stated falling assertion.
\end{proof}

  \medskip
  As an example, consider the cyclic group $\mathbb Z_n$, when $n$
  is odd. The presentation of such a group arises from a word $x^n$ in the free group
  $F_1$ that leads to varieties $V_{-1}$ and $V_{1}$
  satisfying the Minus Condition, and as a consequence we
  have that $P_2(R_A(\mathbb Z_n))=(\frac {n-1} {2},0,1)$ for
  $A \in \{SL(2,\mathbb C), PSL(2,\mathbb C)\}$. So the representation variety over
  $PSL(2,\mathbb C)$ of $\mathbb Z_n$ falls relative to the representation variety
  of $\mathbb Z_n$ over $SL(2,\mathbb C)$. To find examples of where the
  representation variety falls
  properly, torus groups like $<x,y \, ; \, x^py^t=1>$, where $p\geq4$ is even
  and $t\geq3$ is odd, may be used: Theorem 4.7 makes the top-dimensional
  component count strictly smaller. In contrast, notice
  that the representation variety of $\mathbb Z_2$ over $PSL(2,\mathbb C)$ does not fall properly
  or otherwise relative to the representation variety of $\mathbb Z_2$ over $SL(2,\mathbb C)$.

\end{document}